\documentclass[12pt]{amsart}\newtheorem{thm}{Theorem}[section]
\usepackage{amscd}
\usepackage{amssymb}

\newcommand{\Aut}{\operatorname{Aut}}

\newcommand{\bs}{\mathbf{s}}
\newcommand{\bt}{\mathbf{t}}

\newcommand{\bx}{\mathbf{x}}
\newcommand{\cft}{CohFT}
\newcommand{\cfts}{CohFTs}
\newcommand{\CP}[1]{\mathbb{CP}^{#1}}

\newcommand{\F}{\mathcal{F}} 
\newcommand{\G}{\mathcal{G}}    

\newcommand{\HH}{\mathcal{H}}   

\newcommand{\K}{\mathcal{K}}
    
\newcommand{\LL}{\mathcal{L}}   
\newcommand{\mf}{\mathbf{m}}    
\newcommand{\M}{\overline{\MM}} 
\newcommand{\MM}{\mathcal{M}}   
\newcommand{\MP}{\overline{\mathsf{M}}} 
\newcommand{\MPT}{\widetilde{\mathsf{M}}} 
\newcommand{\MT}{\widetilde{\mathcal{M}}} 

\newcommand{\pf}{\mathbf{p}}

\newcommand{\R}{\mathcal{R}}    
\newcommand{\rf}{\mathbf{r}}

\renewcommand{\sf}{\mathbf{s}}
\newcommand{\st}{\operatorname{st}}
\newcommand{\tf}{\mathbf{t}}
\newcommand{\tw}{\tilde{t}}
\newcommand{\xf}{\mathbf{x}}

\newcommand{\del}{\partial} 

\newcommand{\GD}{{\mathbb{G}}} 
\newcommand{\gaf}{\boldsymbol{\gamma}}
\newcommand{\kaf}{\boldsymbol{\kappa}} 
\newcommand{\pit}{\widetilde{\pi}} 
\newcommand{\sif}{\boldsymbol{\sigma}}
\newcommand{\tauf}{\boldsymbol{\tau}}
\newcommand{\zef}{\boldsymbol{0}} 

\newcommand{\nc}{{\mathbb{C}}}  
\newcommand{\nq}{{\mathbb{Q}}}  
\newcommand{\nz}{{\mathbb{Z}}}  

\newcommand{\dis}{\displaystyle}

\newcommand{\la}{\langle}  
\newcommand{\ra}{\rangle}  

\newtheorem{lm}[thm]{Lemma}
\newtheorem{prop}[thm]{Proposition}

\theoremstyle{definition}

\newtheorem{df}[thm]{Definition}

\theoremstyle{remark}

\newtheorem{nota}{Notation} 
\newtheorem{rem}{Remark}  
\newtheorem{ack}{Acknowledgment}  

\begin{document}
\addtocounter{section}{-1}

\title[A Change of Coordinates]
{A Change of Coordinates on the Large Phase Space of Quantum Cohomology} 

\author
[A. Kabanov]{Alexandre Kabanov}
\address
{Mathematik Departement,
ETH-Zentrum,
R\"amistrasse 101, 
CH-8092 Zurich,
Switzerland }
\email{kabanov@math.ethz.ch}

\author
[T. Kimura]{Takashi Kimura}
\address
{Department of Mathematics, 111 Cummington Street, Boston
University, Boston, MA 02215, USA} 
\email{kimura@math.bu.edu}

\thanks{Research of the first author was partially supported by NSF
  grant number DMS-9803553.}

\thanks{Research of the second author was partially supported by NSF grant
  number DMS-9803427.}

\date{\today}

\begin{abstract} 
The Gromov-Witten invariants of a smooth, projective variety $V$, when twisted
by the tautological classes on the moduli space of stable maps, give rise to
a family of cohomological field theories and endow the base of the family
with coordinates. We prove that the potential functions associated to  the
tautological $\psi$ classes (the large phase space) and the $\kappa$ classes
are related by a change of coordinates which generalizes a change of basis
on the ring of symmetric functions.  Our result is a generalization of the
work of Manin--Zograf who studied the case where $V$ is a point. We utilize
this change of variables to derive the topological recursion relations
associated to the $\kappa$ classes from those associated to the $\psi$
classes. 
\end{abstract}

\maketitle

\section{Introduction} 
\label{intro} 
\begin{nota}
All (co)homology are with $\nq$ coefficients unless explicitly
mentioned otherwise. Summation over repeated upper and lower indices
is assumed.   
\end{nota}

\bigskip

The theory of Gromov--Witten invariants of a smooth projective variety $V$ has
developed at a rapid pace \textsl{c.f.}\  \cite{Be,BMa,KM1,RT}. These are
multilinear operations on the cohomology $H^\bullet(V)$ which can be
constructed from intersection numbers on the moduli space of stable maps into
$V$, $\M_{g,n}(V)$. In particular, 
the genus zero Gromov--Witten invariants endow $H^\bullet(V)$ with the
structure of the quantum cohomology ring of $V$. The existence of these
invariants was foreseen by physicists  who encountered these operations as
correlators of a topological sigma model coupled to topological gravity
\cite{W}. These invariants are of great mathematical interest, for example,
because they are symplectic invariants of $V$ \cite{LT} and because of their
close relationship to problems in enumerative geometry \cite{KM1}. 

Gromov--Witten invariants satisfy relations (factorization identities)
parametrized by the relations between cycles on the moduli space of stable
curves $\M_{g,n}$. These relations can be formalized by stating that the space
$(H^\bullet(V),\eta)$ (where $\eta$ is the Poincar\'e pairing) is endowed
with the structure of a cohomological field theory (\cft) in the sense of
Kontsevich--Manin \cite{KM1}. The Gromov--Witten invariants are characterized by 
its generating function (the small phase space potential) $\Phi(\bx)$ where
$\bx\,:=\,\{\,x^\alpha\,\}$ are coordinates associated to a  basis on
$H^\bullet(V)$.  Restricting to genus zero, $\Phi(\bx)$ essentially endows
$(H^\bullet(V),\eta)$ with the structure of a (formal) Frobenius manifold
\cite{Du,Hi,Ma}. It is precisely the structure of a \cft\ which was used by
Kontsevich--Manin to compute the number of rational curves on $\CP{2}$
\cite{KM1} and the number of elliptic curves by Getzler \cite{G} (where the
number is counted with suitable multiplicities).

Furthermore, there are tautological cohomology classes (denoted by $\psi_i$)
associated to the universal curve on $\M_{g,n}(V)$ for all
$i\,=\,1,\,\ldots,\,n$ which are the first Chern class of  tautological line
bundles over $\M_{g,n}(V)$. These classes are a generalization of the $\psi$
classes on $\M_{g,n}$ due to Mumford. What is remarkable is that by
twisting the Gromov--Witten invariants by these $\psi$ classes to obtain the
so-called gravitational descendents, one endows
$(H^\bullet(V),\eta)$ with the structure of a formal family of \cft\
structures whose base is equipped with coordinates
$\bt\,:=\,\{\,t_a^\alpha\,\}$ where $a\,\geq\,1$ and $\alpha$ is as above.
The associated generating function $\F(\bx;\bt)$ (the large phase space
potential) reduces to $\Phi(\bx)$ when $\bt$ vanishes. The large phase space
potential $\F$ is itself a remarkable object as its exponential is
conjectured to satisfy a highest weight condition for the Virasoro algebra
\cite{E}, a conjecture which has nontrivial consequences \cite{GePa}. Indeed,
when $V$ is a point, this condition is equivalent to the Witten conjecture
\cite{W} proven by Kontsevich \cite{Ko}.

There are other tautological cohomology classes on $\M_{g,n}(V)$ associated
to its universal curve. In this paper, we define generalizations to
$\M_{g,n}(V)$ of the ``modified'' $\kappa$ classes (due to Arbarello--Cornalba
\cite{AC}) on $\M_{g,n}$.  We define Gromov--Witten invariants twisted by the
$\kappa$ classes and prove that we obtain a formal family of \cfts\ on
$(H^\bullet(V),\eta)$ whose base is endowed with coordinates $\bs\,:=\,
\{\,s_a^\alpha\,\}$ where $a\,\geq\,0$. We then prove that the generating
function $\G(\bx;\bs)$ associated to this family can be identified 
with the large phase space potential $\F(\bx;\bs)$ through an explicit change
of variables.  This change of variables can be interpreted as a change of
basis in the space of symmetric functions whose variables take values in
$H^\bullet(V)$.  The variables $\bs$ can be interpreted as another canonical
set of coordinates on the large phase space. We also utilize this change of
variables to derive topological recursion relations for  $\G$ in terms of
those of $\F$. 

When $V$ is convex, the $\kappa$ classes on $\M_{0,n}(V)$ had already been 
introduced in \cite{KK2} and the genus zero topological recursion relations
were proven. This paper generalizes those results to situations where
$\M_{g,n}(V)$ need not have the expected dimension (and, hence, the 
technicalities of the virtual fundamental class cannot be avoided) as
well as to derive the change of variables on the large phase space.

When $V$ is a point, our formula reduces to the work of
Kaufmann--Manin--Zagier \cite{KMZ} and \cite{MaZo} who noted (see also
\cite{KMK}) that, in addition, the coordinates $\bs$ are additive with
respect to the tensor product in the category of \cfts. Manin--Zograf
\cite{MaZo} used this formula to compute asymptotic Weil--Peterson
volumes of the moduli spaces $\M_{g,n}$ as $n\,\to\,\infty$ (this was done
for $g\,=\,1$ in \cite{KK}). However, this additivity property need not hold
for a general variety $V$.

It is worth pointing out several generalizations. First of all, when $V$ is a
point, Manin--Zograf use the Witten conjecture to show that their change of
variables can be directly interpreted as arising from an analogous change of
the cohomology classes appearing in the potential functions. It would be
interesting to obtain an analogous result for a general $V$. Secondly, the
above construction should be feasible for any \cft\ and there should be
coordinates which are additive under tensor product -- such a construction
would be useful in studying the ring of \cfts. Work towards this direction is
in progress \cite{KK3}. The third is the fact that there are yet another set of
tautological classes (called $\lambda$) on $\M_{g,n}$ associated to the Hodge
bundles. Twisting Gromov--Witten invariants by both the $\kappa$ and $\lambda$
classes, one obtains the very large space \cite{KK,MaZo} (see also
\cite{FaPa}),  a subset of which form coordinates on the moduli space of
nondegenerate rank one \cfts\  in genus $1$ \cite{KK} which are additive under
tensor product. It would be interesting to understand the role of these
additional coordinates for general $V$.

The first section of the paper is a review of the technicalities necessary to 
push forward and pull back cohomology classes on the moduli space of stable
maps. This includes Gysin morphisms and the flat push-forward.

In the second section, we review the basic properties of the moduli
space of stable maps, the structure of the boundary classes, and
properties of the virtual fundamental classes.

In the third section we introduce the tautological $\kappa$ and $\psi$
classes, and prove its restriction properties on the boundary classes.

In the fourth section, we define the notion of a \cft\ and its potential
function.  We review the large phase space potential $\F$. We prove that by
introducing the $\kappa$ classes, $(H^\bullet(V),\eta)$ is endowed with a
formal family of \cft\ structures together with coordinates on the base of
the family.

In the fifth section, we prove that after an explicit change of
coordinates, the potential $\G$ can be identified with the large phase
space potential $\F$.

In the final section, we derive the topological recursion relations for $\G$
in genus 0 and 1 and derive the usual topological recursion relations for
$\F$ through the change of variables.

\begin{sloppypar}
\begin{ack}
We would like to thank D.~Abramovich for useful conversations. 
\end{ack}
\end{sloppypar}

\section{Technical Preliminaries.}
\label{tech}

In this section we present several technical points needed in the
sequel. They are concerned with the Gysin morphisms in homology and
cohomology. You may skip this section provided you are willing to
accept that everything works at a ``naive'' level. An article of
Fulton and MacPherson \cite{FuMP} may serve as a general reference to
this section. All references mentioned in this section deal with
schemes rather than stacks, but the sheaf-theoretic approach allows
one to work in the category of stacks. If $F$ is a functor, then $RF$
denotes the corresponding derived functor.

Let $\pi: Y \to X$ be a flat representable morphism of
Deligne--Mumford stacks with fibers of pure dimension $d$. As explained
in \cite{De} $\pi$ defines the natural morphism of $Tr_\pi: R^{2d}
\nq \to \nq$, which induces the corresponding flat push-forward in
cohomology with compact supports $\pi_*: H_c^k (Y) \to H_c^{k-2d}
(X)$. (The axioms uniquely defining the morphism $Tr$ are also given
in \cite{Ve1} and \cite{La}.)

One of the axioms defining the $Tr$ morphism states that it commutes
with the base change, that is, with the pull-back on cohomology in a
fibered square. However, in this paper we will need to consider
commutative squares which are a little more general than the fibered
squares. It is the reason for giving the following definition.

\begin{df}\label{df:fibered}
  Let $X_1,Y_1,X,Y$ be Deligne--Mumford stacks. A commutative square 
\[
\begin{CD}
  Y_1 @>f_1>> Y \\
  @VV\pi_1 V @VV\pi V \\
  X_1 @>f>>  X \\
\end{CD}
\]
is called \textit{close to a fibered square} if the induced morphism
$g: Y_1 \to X_1 \times_X Y$ is a proper birational morphism, and there
is an open subset $U$ of $X_1 \times_X Y$ whose intersection with each
fiber of $\pi_1$ is a dense subset of the fiber such that $g|_{g^{-1}
  U}$ is an isomorphism.
\end{df}

\begin{lm}\label{lm:fibered}
  If the commutative square
\[
\begin{CD}
  Y_1 @>f_1>> Y \\
  @VV\pi_1 V @VV\pi V \\
  X_1 @>f>>  X \\
\end{CD}
\]
is close to a fibered one, and $\pi$ and $\pi_1$ are representable
flat morphisms with fibers of pure dimension $d$, then $\pi_* f^* =
f_1^* \pi_{1*}: H^\bullet (Y) \to H^\bullet (X_1)$.
\end{lm}

\begin{proof}
Consider the following diagram
\[
\begin{CD}
  Y_1 @>g>> X_1 \times_X Y @>pr_2>> Y \\
  @V\pi_1VV  @Vpr_1 VV  @V\pi VV \\
  X_1 @=    X_1            @>f>>    X, \\
\end{CD}
\]
where $pr_2\ g = f_1$. Since the right square is a fibered square it
follows from the properties of the $Tr$ morphism that $\pi_* f^* =
pr_2^*\ pr_{1*}$. Therefore it remains to show that $pr_{1*} = g^*
\pi_{1*}$. It follows from the construction in \cite[Sec.~2]{De} that
the $Tr$ morphism is determined by a Zariski open subset whose
intersection with each fiber is dense. In other words, $Tr: R\ pr_{1!}
\ \nq \to \nq$ coincides with 
\[
R\ pr_{1!}\ \nq \to R\ pr_{1!}\ Rg_*\ \nq = R \pi_{1!}\ \nq \to \nq.
\]
We have used the fact that $R\ g_* =
R\ g_!$ since $g$ is proper.
\end{proof}

Dually, a flat morphism $\pi: Y \to X$ with fibers of pure dimension
$d$ determines a flat pull-back $\pi^*: H_k (X) \to H_{k+2d} (Y)$.  It
is shown in \cite[Sec.~6]{La} that $\pi^*$ agrees with the flat pull-back
$\pi^*: A_k (X) \to A_{k+d} (Y)$ via the cycle map. (In the set up of
bivariant intersection theory \cite[2.3]{FuMP} each flat morphism
$\pi$ determines a canonical element in $T^{-2d} (Y \to X)$.)

We also need to define the Gysin morphisms associated to regular
imbeddings. A closed imbedding $i:X_1 \to X$ is called a
\textit{regular imbedding of codimension $d$} if the conormal sheaf of
$X_1$ in $X$ is a locally free sheaf on $X_1$ of rank $d$
\cite[B.7.1]{Fu}. Let $i: X_1 \to X$ be a regular embedding of
codimension $d$. The corresponding canonical element $\theta_i \in
H^{2d} (X, X-X_1)$ is constructed in \cite[IV.4]{BFM} and
\cite[Sec.~5]{Ve2}. (In bivariant intersection theory $H^{2d}(X,X-X_1)
= T^{2d}(X_1 \to X)$.) If 
\[
\begin{CD}
  Y_1 @>i_1>> Y \\
  @Vf_1VV  @VfVV \\
  X_1 @>i>>  X\\
\end{CD}
\]
is a fibered square, then the pull-back $f^*\theta_i$ determines an
element in $H^{2d} (Y,Y-Y_1)$. Accordingly, it defines Gysin
homomorphisms: 
\begin{align*}
  & i^!: H_k(Y) \to H_{k-2d}(Y_1) \quad \text{and} \\
  & i_!: H^k(Y_1) \to H^{k+2d}(Y) 
\end{align*}
by the cap-product (or cup-product) with $f^* \theta_i$.  However, we
will denote $i^!: H_k(X) \to H_{k-2d}(X_1)$ by $i^*$, and $i_!:
H^k(X_1) \to H^{k+2d}(X)$ by $i_*$. This agrees with the notation from
\cite{Fu}. If $E\to X$ is a rank $d$ vector bundle, and $X_1$ is the
zero scheme of a section $i:X \to E$, then $i_* 1 = c_d E$
\cite[Sec.~19.2]{Fu}. The Gysin morphism $i^!$ defined above agrees
with the Gysin morphism $i^!$ on the level of Chow groups via the
cycle map \cite{Ve2}.

\begin{rem}
  More generally, one can define the Gysin morphisms for local
  complete intersection morphisms. If a flat morphisms is at the same
  time a local complete interesection morphism, then two definitions
  agree. 
\end{rem}

The morphisms $\pi_*,\pi^*,i_!,i^!$ satisfy the expected projection
type formulae and commute with the standard pull-backs and
push-forwards \cite[2.5]{FuMP}. We will use these properties without
explicitly mentioning them.

\section{The Moduli Spaces of Stable Maps}
\label{define}

\bigskip

We adopt the notation from \cite{Ge}. Let $\M_{g,n}$ be the moduli
space of stable curves. The stability implies that $2g-2+n>0$. Let
$\Gamma$ be a stable graph of genus $g$ with $n$ tails. We denote by
$\M(\Gamma) \subset \M_{g,n}$ the closure in $\M_{g,n}$ of the locus
of stable curves with the dual graph $\Gamma$, and by $i_\Gamma$ the
corresponding inclusion. Let
\[
\MT(\Gamma) := \prod_{v\in V(\Gamma)} \M_{g(v),n(v)}. 
\]
Then $\Aut(\Gamma)$ acts on $\MT(\Gamma)$. The natural morphism
\[
\mu_\Gamma: \MT(\Gamma) \to \M(\Gamma)
\]
identifies $\M(\Gamma)$ with $\MT(\Gamma) / \Aut(\Gamma)$. We denote
by $\rho_\Gamma$ the composition
\begin{equation*}
  \rho_\Gamma: \MT(\Gamma) \stackrel{\mu_\Gamma}{\longrightarrow}
  \M(\Gamma) \stackrel{i_\Gamma}{\longrightarrow} \M_{g,n}. 
\end{equation*}

The previous considerations apply word for word to the moduli spaces
of prestable curves $\MP_{g,n}$, $g\ge 0, n\ge 0$, their subspaces
$\MP(\Gamma)$, and the products $\MPT(\Gamma)$ \cite[Sec.~2]{Ge}. Note
that $\M_{g,n}$ is an open dense substack of $\MP_{g,n}$ when
$2g-2+n>0$, and, more generally, $\M(\Gamma)$ is an open dense
substack of $\MP(\Gamma)$ when $\Gamma$ is a stable graph.

We adopt a similar notation for the substacks of $\M_{g,n}(V,\beta)$
determined by decorated stable graphs. Let $H_2^+(V,\nz)$ denote the
semigroup generated by those homology classes represented by the image
of a morphism from a curve into $V$. Let $\GD$ be a stable graph of
genus $g$ with $n$ tails whose vertices are decorated by elements
of $H_2^+(V,\nz)$.  (Henceforth, such decorated graphs will be denoted by
$\GD$.) Then we denote by $\M(\GD,V)$ the closure in $\M_{g,n}(V,\beta)$ of
those points in the moduli space of stable maps whose dual graph is $\GD$.
Let $\MT(\GD,V)$ be determined by the following fibered square
(cf.~\cite[Sec.~6]{Ge}):
\[
\begin{CD}
  \MT(\GD,V) @>\Delta>> \dis{\prod_{v\in V(\GD)} }
  \M_{g(v),n(v)}  (V,\beta(v)) \\
  @VVev V @VVev V\\
  V^{E(\GD)} @>\Delta_1>> V^{E(\GD)} \times V^{E(\GD)},
\end{CD}
\]
where $\Delta_1$ is the diagonal morphism. In the sequence
\begin{equation*}
  \prod_{v\in V(\GD)} \M_{g(v),n(v)} (V,\beta(v)) 
  \stackrel{\Delta}{\longleftarrow} 
  \MT(\GD,V) \stackrel{\mu(\GD)}{\longrightarrow} \M(\GD,V)
  \stackrel{i(\GD)}{\longrightarrow} \M_{g,n} (V,\beta),
\end{equation*}
the morphism $\mu(\GD)$ is the quotient by $\Aut(\GD)$ identifying
$\M(\GD,V)$ with the quotient, and $i(\GD)$ is the inclusion of a
substack. We denote the composition of two morphisms on the right by
$\rho(\GD)$.

We will also need to introduce some other notation to describe the pull back
of the 
virtual fundamental classes with respect to the inclusions of the
strata (cf.~\cite[Sec.~6]{Ge}. Let $\Gamma$ be a graph of genus $g$
with $n$ tails, not necessarily stable. We define
\[
\M(\Gamma,V,\beta) :=
 \M_{g,n}(V,\beta) \times_{\MP_{g,n}} \MP(\Gamma). 
\]
It is the closure of the subset of $\M_{g,n}(V,\beta)$ whose points
correspond to the graph $\Gamma$ after forgetting the decoration. If
$\Gamma$ is a stable graph, then $\M(\Gamma,V,\beta) =
\M_{g,n}(V,\beta) \times_{\M_{g,n}} \M(\Gamma)$ since $\M(\Gamma)$ is
dense in $\MP(\Gamma)$. If $\GD$ is a decorated graph we denote by
$\GD^0$ the underlying non-decorated graph. Let
\[
\MT(\Gamma,V,\beta) := \coprod_{\GD: \GD^{0} = \Gamma} 
\MT(\GD,V).  
\]
As before, $\Delta: \MT(\Gamma,V,\beta) \to
\coprod_{\GD^{0}=\Gamma} \prod_{v\in V(\GD)} \M_{g(v),n(v)}
(V,\beta(v))$ is determined by the diagonal morphism. One has the
natural morphism:
\[
\rho(\Gamma): \MT(\Gamma,V,\beta) 
\stackrel{\mu(\Gamma)}{\longrightarrow} \M(\Gamma,V,\beta)
\stackrel{i(\Gamma)}{\longrightarrow} \M_{g,n} (V,\beta). 
\]
Here $i(\Gamma)$ is an inclusion of a substack, and $\mu(\Gamma)$
factors as 
\[
\MT(\Gamma,V,\beta) \to \MT(\Gamma,V,\beta) / \Aut(\Gamma) \to
\M(\Gamma,V,\beta), 
\]
where the second morphism is a proper, surjective, birational
morphism.  The difference with the previous situation is explained by
the fact that two substacks $\M(\GD,\beta)$ and $\M(\GD',\beta)$ of
$\M_{g,n}(V,\beta)$ whose underlying undecorated graphs are the same
may have a nonempty intersection.

\section{Tautological Classes}
\label{taut}

In this section we introduce the tautological $\kappa$ classes on the
moduli spaces of stable maps which generalize the corresponding
tautological classes on the moduli spaces of stable curves. We will
also show how these classes restrict to the boundary strata.

Let $\pi: \M_{g,n+1}(V, \beta) \to \M_{g,n} (V,\beta)$ be the
universal curve. We assume that $\pi$ ``forgets'' the
$(n+1)^{\text{st}}$ marked point. The morphism $\pi$ has $n$ canonical
sections $\sigma_1, \ldots, \sigma_n$. Each of these sections
determines a regular embedding. We denote by $\omega$ the relative
dualizing sheaf of $\pi$.

\begin{df}
  For each $i=1,\ldots, n$ the {\it tautological line bundle} $\LL_i$
  on $\M_{g,n}(V,\beta)$ is $\sigma_i^* \omega$. The {\it tautological
    class} $\psi_i \in H^2(\M_{g,n}(V,\beta))$ is the first Chern
  class $c_1(\LL_i)$.
\end{df}

\begin{rem}
  It is shown in \cite[Sec.~5]{Ge} that $\psi_i = p^* \Psi_i$, where
  $p: \M_{g,n} (V,\beta) \to \MP_{g,n}$, where $\Psi_i$ is the
  tautological class in $H^2(\MP_{g,n})$.  
\end{rem}

One can also pull back cohomology classes from $V$ to $\M_{g,n}(V,\beta)$
using the evaluation maps to obtain the Gromov--Witten classes.
The definition of the $\kappa$ classes involves both, powers of the $\psi$
classes and these pull backs. 

\begin{df}
  The {\it tautological class} $\kappa_{a}$ in
  $H^{\bullet}(\M_{g,n}(V,\beta))\, \otimes\,H^\bullet(V)^*$ for $a\ge
  -1$ is defined as follows. For each $\gamma\,\in\,H^\bullet(V)$, the
  cohomology class $\kappa_{a}(\gamma)$ is $\pi_* (\psi_{n+1}^{a+1}
  ev_{n+1}^* (\gamma) )$, where $\pi$ is the universal curve defined
  above. In particular, if $\gamma$ has definite degree $|\gamma|$
  then $\kappa_{a} (\gamma)$ has degree $2a + |\gamma|$. If $\{\,
  e_\alpha\, \}_{\alpha \in A}$, is a homogeneous basis for
  $H^\bullet(V)$, then $\kappa_{a,\alpha}$ denotes the cohomology
  class $\kappa_a (e_\alpha)$.
\end{df}

\begin{rem}
  The class $\kappa_{-1}(\gamma)$ vanishes due to dimensional reasons
  if $|\gamma|<2$. In addition, all classes $\kappa_{-1}(\gamma)$
  vanish on $\M_{g,n}(V,0)$.  The classes $\kappa_{-1} (\gamma)$ are
  not needed in the change of coordinates formula in
  Sec.~\ref{change}.
\end{rem}

Our definition corresponds to the ``modified'' $\kappa$ classes
defined by Arbarello and Cornalba \cite{AC} rather than the
``classical'' $\kappa$ classes defined by Mumford \cite{Mu}.

The following lemma shows how the $\kappa$ classes restrict to the
boundary substacks of $\M_{g,n} (V,\beta)$. 

\begin{lm}
\label{kappa:restr}
Let $\GD$ be a stable $H_2^+(V,\nz)$ decorated genus $g$, degree $\beta$
graph with $n$ tails. Denote the class $\kappa_{a} (\gamma)$ on
$\M_{g,n} (V,\beta)$ (resp. $\M(v)$, where $v\in V(\GD))$ by $\kappa$
(resp. $\kappa_v)$. Then
\[
\rho(\GD)^* (\kappa) = \Delta^* \sum_{v\in V(\GD)} \kappa_v.
\]
\end{lm}

\begin{proof}
  Let $v\in V(\GD)$. Denote by $\GD(v)$ the graph obtained from $\GD$
  by attaching a tail labeled $n+1$ to the vertex $v$ of $\GD$.  For
  each $v\in V(\GD)$ the graph $\GD(v)$ determines a substack of
  $\M_{g,n+1} (V,\beta)$, and there are natural morphisms
\[
\pit: \prod_{w\in V(\GD(v))} \M(w) \to \prod_{w\in V(\GD)} \M(w), 
\ \text{and} \
\pit: \MT(\GD(v),V) \to \MT(\GD,V). 
\]

Consider the following commutative diagram
\[
\begin{CD}
\dis{\coprod_{v\in V(\GD)}} \prod_{w\in \GD(v)} \M(w) @<\coprod\Delta<<
\dis{\coprod_{v\in V(\GD)}} \MT(\GD(v),V) @>\coprod\rho(\GD')>>
\M_{g,n+1} (V,\beta) \\
@VV \coprod\pit V @VV \coprod\pit V @VV\pi V \\
\prod_{w\in \GD} \M(w) @< \Delta<<
\MT(\GD,V) @>\rho(\GD)>>
\M_{g,n} (V,\beta). \\
\end{CD}
\]
Note that the left square is a fibered square, the right square is
close to a fibered square in the sense of Def.~\ref{df:fibered}, and
all morphisms $\pi$, $\pit$ are representable and flat. Therefore, one
can apply Lem.~\ref{lm:fibered}. Also note that for each $v\in V(\GD)$
one has
\[
\rho(\GD(v))^* \psi_{n+1} = \Delta^* \psi_{n+1} \quad \text{and} 
\quad \rho(\GD(v))^* ev_{n+1}^* \gamma = \Delta^* ev_{n+1}^* \gamma. 
\]
Now
\begin{equation*}
\begin{split}
&\rho(\GD)^* (\kappa) =
\rho(\GD)^* \pi_* (\psi_{n+1}^{a+1} ev_{n+1}^* \gamma ) = 
\sum_{v\in V(\GD)} \pit_* \rho(\GD(v))^* (\psi_{n+1}^{a+1} ev_{n+1}^*
      \gamma ) \\ 
&= \sum_{v\in V(\GD)} \pit_* \Delta^* (\psi_{n+1}^{a+1} ev_{n+1}^*
       \gamma ) =  
\sum_{v\in V(\GD)} \Delta^* \pit_* (\psi_{n+1}^{a+1} ev_{n+1}^* \gamma)=  
\Delta^* \sum_{v\in V(\GD)} \kappa_v. 
\end{split}
\end{equation*}
\end{proof}

The above lemma shows that the class $\kappa_a(\gamma)$ restricts to
the sum of the $\kappa_a(\gamma)$ classes. It follows that
$\exp(\kappa_a(\gamma))$ restricts to the product of
$\exp(\kappa_a(\gamma))$. More generally, $\exp(\sum_{a=-1}^\infty
\kappa_{a,\alpha} s_a^\alpha)$, where $s_i$'s are formal variables,
restricts to the product of $\exp(\sum_{a=-1}^\infty \kappa_{a,\alpha}
s_a^\alpha)$. This will be used in Sec.~\ref{cohft}.

\section{Cohomological Field Theories}
\label{cohft}

In this section, we define a cohomological field theory in the sense
of Kontsevich--Manin \cite{KM1}. We prove that the Gromov--Witten
invariants twisted by the $\kappa$ classes endows $H^\bullet(V)$
together with its Poincar\'e pairing with the a family of \cft\ 
structures. In genus zero, this reduces to endowing $H^\bullet(V)$
with a family of formal Frobenius manifold structures arising from the
Poincar\'e pairing and deformations of the cup product on
$H^\bullet(V)$. These deformations contain quantum cohomology as a
special case.

\begin{df}  Let $(\HH,\eta)$ be an $r$-dimensional vector space $\HH$ with an
  even, symmetric nondegenerate,  bilinear form $\eta$. A \textsl{(complete)
 rank $r$ cohomological field theory (or \cft) with state space $(\HH,\eta)$}
 is a collection  
$\Omega\,:=\,\{\,\Omega_{g,n}\,\}$ where $\Omega_{g,n}$ is an even element in
$\R_{g,n}\,:=\,H^\bullet(\M_{g,n}) \otimes T^n\HH^*$ defined for stable pairs
$(g,n)$ satisfying (i) to (iii) below (where the summation convention has been
used): 
\begin{description}
\item[i] $\Omega_{g,n}$ is invariant under the diagonal action of the
  symmetric   group $S_n$ on $T^n\,\HH$ and $\M_{g,n}$.
\item[ii] For each partition of $[n]\,=\,J_1\sqcup J_2$ such that $|J_1|=n_1$
 and  $|J_2|=n_2$ and nonnegative $g_1$, $g_2$ such that $g\,=\,g_1\,+\,g_2$
 and $2g_i-2+n_i+1\,>\,0$ for all $i$,  consider the inclusion map 
$\rho\,:\,\M_{g_1,J_1 \sqcup *}\times\M_{g_2,J_2\sqcup
 *}\,\to\,\M_{g_1+g_2,n}$ where $*$ denotes the two marked points that are
 attached under the inclusion map. The forms satisfy the restriction property
\begin{align*}
&\rho^*\,\Omega_{g,n}(\gamma_1,\gamma_2\,\ldots,\gamma_n)\,=\, \\
&\pm\,\Omega_{g_1,n_1}((\bigotimes_{\alpha\in J_1}\,\gamma_\alpha)\,\otimes
e_\mu)\,\eta^{\mu
  \nu}\,\otimes\,\Omega_{g_2,n_2}(e_\nu\,\otimes\,\bigotimes_{\alpha\in 
  J_2}\,\gamma_\alpha) 
\end{align*}
where the sign $\pm$ is the usual one obtained by applying the permutation
induced by the partition to $(\gamma_1,\gamma_2,\ldots\,\gamma_n)$
taking into account the grading of $\{\,\gamma_i\,\}$ and where
$\{\,e_\alpha\,\}$ is a homogeneous basis for $\HH$.

\item[iii] Let $\rho_0\,:\,\M_{g-1,n+2}\,\to\,\M_{g,n}$ be the canonical map
  corresponding to attaching the last two marked points together then
\begin{equation*}
\rho_0^*\,\Omega_{g,n}(\gamma_1,\gamma_2,\ldots,\gamma_n)\,=\,
\Omega_{g-1,n+2}\,(\gamma_1,\gamma_2,\ldots,\gamma_n, e_\mu,
e_\nu)\,\eta^{\mu\nu}. 
\end{equation*}

\item[iv] If, in addition, there exists an even element $e_0$ in $\HH$
such that 
\[
\pi^*\,\Omega_{g,n}(\gamma_1,\ldots,\gamma_n)\,=\,\Omega_{g,n+1}(\gamma_1,
\ldots, \gamma_n, e_0)
\]
and
\[
\int_{\M_{0,3}}\,\Omega_{0,3}(e_0,\gamma_1,\gamma_2)\,=\,
\eta(\gamma_1,\gamma_2)
\]
for all $\gamma_i$ in $\HH$ then $\Omega$ endows $(\HH,\eta)$ with the
  structure of a \textsl{\cft\ with flat identity $e_0$.}
\end{description}

A \emph{cohomological field theory of genus $g$} consists of only those
$\Omega_{g',n}$ where $g'\,\leq\,g$ which satisfy the subset of axioms of a
cohomological field theory which includes only objects of genus $g' \leq g$.
\end{df}

The strata maps $\rho$ and $\rho_0$ in the above definition can be
extended to arbitrary boundary strata on $\M_{g,n}$. Let $\Gamma$ be a stable
graph then there is a canonical map $\rho_\Gamma$  obtained by composition of
the canonical maps 
\[
\prod_{v\in V(\Gamma)}\,\M_{g(v),n(v)}\,\to\,
\M_\Gamma\,\to\,\M_{g,n}.
\]
Since the map $\rho_\Gamma$ can be constructed from morphisms in $(ii)$ and
$(iii)$ above, $\Omega_{g,n}$ satisfies a restriction property of the form
\begin{equation} \label{restriction}
\rho_\Gamma^*\Omega_{g,n}\,=\, \eta^{-1}_\Gamma(\,\bigotimes_{v\in V(\Gamma)}
\Omega_{g(v),n(v)}\,)
\end{equation}
where 
\[ \eta^{-1}_\Gamma\,:\,\bigotimes_{v\in V(\Gamma)} \R_{g(v),n(v)}\, \to\,
\R_{g,n} \]
is the linear map contracting tensor factors of $\HH$ using the metric $\eta$
induced from successive application of equations $(ii)$ and $(iii)$ above.

Notice that the definition of a cohomological field theory is valid even when
enlarges the ground ring from $\nc$ to another ring $\mathcal{K}$. 

Finally, axioms $(i)$ to $(iii)$ in the definition of a \cft\ is equivalent
to endowing $(\HH,\eta)$ with the structure of an algebra over the modular
operad $H_\bullet(\M)\,:=\,\{\,H_\bullet(\M_{g,n})\,\}$ \cite{GK}. 


\begin{df}
Let $\Lambda$ consist of formal symbols $q^\beta$ for all
$\beta\,\in\,H_2^+(V,\nz)$ together with the multiplication
$(q^\beta\,q^{\beta'})\,\mapsto\, q^{\beta+\beta'}$. Let
$\nc[[\Lambda]]$ consist of formal sums
$\sum_{\beta\,\in\, H_2^+(V,\nz)}\,a_\beta\,q^\beta$ where $a_\beta$ are
elements in $\nc$. Assign to each $q^\beta$, the degree $-2
c_1(V)\,\cap\,\beta$. The product is well-defined according to
\cite[Prop.~II.4.8]{Kol}. This endows $\Lambda$ with the structure of
a semigroup with unit. Furthermore, let $\nc[[\Lambda,\sf]]\, := \,
\nc[[\Lambda]][[\sf]]$, formal power series in the variables $\sf$
with coefficients in $\nc[[\Lambda]]$.
\end{df}

\begin{nota}
Let $V$ be a topological space and let $H^\bullet(V,\nc)$ be given a
homogeneous basis $\mathbf{e}\,:=\,\{\,e_\alpha\,\}_{\alpha\in A}$ and let
$e_0$ denote the identity element. Let
$\sf\,:=\,\{\,s_a^\alpha\,|\,a\geq -1, \alpha\in A\,\}$ be a
collection of formal variables with grading $|s_a^\alpha|
\,=\,2a+|e_\alpha|$. All formal power series and polynomials in a collection
of variables (\emph{e.g.}\ $\sf$) are in the $\nz_2$-graded sense.
\end{nota}

It will be useful to associate a generating function (called the
\textsl{potential}) to each \cft.

\begin{df}
Let $\Omega$ be a rank $r$ \cft\ with state space $(\HH,\eta)$. Its
  \emph{potential function} $\Phi$ in $\lambda^{-2}\,\nc[[\HH,\lambda]]$ is
  defined by
\[
\Phi(\bx)\,:=\,\sum_{g\,\geq\,0}\,\Phi_g(\bx)\,\lambda^{2g-2}
\]
where 
\[
\Phi(\bx) := \sum_{n=3}^\infty
\, \frac{1}{n!}\,\int_{\M_{g,n}}
\Omega_{g,n}(\bx,\bx,\ldots,\bx)
\]
and $\bx\,=\,\sum_{\alpha\,=\,0}^{r-1}\, x^\alpha\,e_\alpha$ for a given
homogeneous basis $\{\,e_0,\ldots,e_{r-1}\,\}$ for $\HH$. The formal parameter
$\lambda$ is even.
\end{df}

In genus zero, the potential function yields yet another formulation of a
\cft\  which is essentially the definition of a formal Frobenius manifold
structure on its state space. 

\begin{thm}
\label{thm:wdvvg}
Let $(\HH,\eta)$ be an $r$ dimensional vector space with metric. An element
$\Phi_0(\bx)$ in $\nc[[\HH]]$ is the potential of a rank $r$, genus zero
\cft\ $(\HH,\eta)$ if and only if \cite{KM1,Ma} it contains only terms which
are cubic and higher order in the variables $x^0,\ldots,x^r$ and it satisfies
the WDVV equation
\[
(\partial_{a} \partial_{b} \partial_{e}\Phi_0)\, \eta^{ef}\, (\partial_{f}
\partial_{c} \partial_{d}\Phi_0)\, = \,(-1)^{|x_a|(|x_b| + |x_c|)}\,
(\partial_{b} \partial_{c}
\partial_{e} \Phi_0)\, \eta^{ef} \,(\partial_{f} \partial_{a}
\partial_{d}\Phi_0),
\]
where $\eta_{ab} := \eta(e_a,e_b)$, $\eta^{ab}$ is in inverse matrix to
$\eta_{ab}$,  $\partial_a$ is derivative with respect to $x^a$, and the
summation convention has been used. Furthermore, any genus zero \cft\ is
completely characterized by its genus zero potential $\Phi_0(\bx)$.
\end{thm} 

The theorem follows from the work of Keel \cite{Ke} who proved that all
relations between boundary divisors on $\M_{0,n}$ arise from lifting the
basic codimension one relation on $\M_{0,4}$.

As before, one can extend the ground ring $\nc$ above to $\nc[[\Lambda,\sf]]$
in the definition of the potential of a genus zero \cft\  and the above
theorem extends, as well. In our setting, the potential is a formal function
on $\HH\,:=\,H^\bullet(V,\nc[[\Lambda,\sf]])$ and $\eta$ is the Poincar\'e
pairing extended linearly to $\nc[[\Lambda,\sf]]$. $\Phi$ belongs to
$\lambda^{-2}\, \nc[[\Lambda,\sf,\lambda]][[x^0,\ldots,x^r]]$  Again, if
$H^\bullet(V)$ consists entirely of even dimensional classes then plugging in
numbers (almost all of which are zero) for all $s_a^\alpha$ where
$a=-1,0,1,\ldots$ and $\alpha\,=\,0,1,\ldots,r-1$ and setting $\lambda\,=\,1$,
one obtains families of \cft\ structures on $H^\bullet(V,\nc[[\Lambda]])$.

\begin{nota}
  We define $\sf\kaf$ to be $\dis{\sum_{a=-1}^\infty
    \kappa_{a,\alpha}} s_a^\alpha$. Note that each term has even
  parity.
\end{nota}

\begin{thm} \label{thm:cohft}
  Let $V$ be a smooth projective variety. For each pair $(g,n)$ such
  that $2g-2+n>0$, let $\Omega_{g,n}$ be the element of
  $\R_{g,n}(V)[[\Lambda,\sf]]$ defined by 
\begin{equation*} \label{omega}
\begin{split}
\Omega_{g,n} & (\gamma_1,\ldots,\gamma_n)\, \\
& :=\,\sum_{\beta\in
  H_2^+(V,\nz)}\,\st_*(\,ev_1^* 
\gamma_1\,\cdots\, ev_n^*\gamma_n\,\exp(\kaf \sf) \cap
[\M_{g,n}(V,\beta)]^{virt})
\,q^\beta, 
\end{split} 
\end{equation*}
where $\gamma_1, \gamma_2, \ldots, \gamma_n$ are elements in
$H^\bullet(V,\nc)$. Then $\Omega\,:=\,\{\,\Omega_{g,n}\,\}$ endows
$(H^\bullet(V,\nc[[\Lambda,\sf]]),\eta)$ with the structure of a \cft\ 
where $\eta$ is the Poincar\'e pairing extended
$\nc[[\Lambda,\sf]]$-linearly.
\end{thm}

\begin{proof}
  It is clear that the morphisms $\Omega_{g,n}$ are $S_n$-equivariant.
  In order to prove the restriction properties fix $\beta\in
  H_2^+(V,\nz)$, $(g,n)$ such that $2g-2+n>0$, and a stable graph
  $\Gamma$ of genus $g$ with $n$ tails. Let $G$ be the set of all
  $H_2^+(V,\nz)$ decorated graphs such that the underlying graph
  without decoration is $\Gamma$. Let
\[
X := \coprod_{\GD \in G}
X_\GD, \quad \text{where} \quad X_\GD:= \prod_{v\in \GD} \M(v),  
\]
and let $[X_\GD]^{virt}\in H_\bullet (X_\GD)$ be the product of the
corresponding virtual fundamental classes.  Consider the following
commutative diagram:
\begin{equation*}
  \begin{CD}
    X @<\Delta<< \MT(\Gamma,V,\beta) @>\mu(\Gamma)>> 
         \M(\Gamma,V,\beta) @> i(\Gamma)>> \M_{g,n}(V,\beta) \\
    @VV\st V @VV\st V @VV\st V @VV\st V \\
    \MT(\Gamma) @= \MT(\Gamma) @>\mu_\Gamma>> \M(\Gamma) 
         @>i_\Gamma>> \M_{g,n}.
  \end{CD}
\end{equation*}
We want to see how the $\beta$ summand of $\Omega_{g,n}$ restricts to
$H^\bullet(\MT(\Gamma))$. 

In the sequence of equations below we will use the following
properties. The right square of the above diagram is a fibered square.
All vertical morphisms $\st$ are proper. If $x \in H^\bullet
(\MT(\Gamma))$ is invariant under the action of $\Aut(\Gamma)$, then
$\mu_\Gamma^* \mu_{\Gamma *} x = N x$, where $N:=|\Aut(\Gamma)|$. In
addition, we use the following result of Getzler \cite[Thm.~13]{Ge}:
\[
i_\Gamma^! [\M_{g,n}(V,\beta)]^{virt} = 
\frac{1}{N} \mu(\Gamma)_* \Delta_1^! 
\sum_{\GD\in G} [X_\GD]^{virt},  
\]
where 
\[
\Delta_1: V^{E(\GD)} \to  V^{E(\GD)} \times V^{E(\GD)} 
\]
is the diagonal morphism. 

If $\GD \in G$, and $v\in V(\GD)$, then we denote by $\gaf_v$ the
tensor product of the corresponding $\gamma_i$'s on $\M(v)$, and by
$\kaf_v \sf$ the formal sum on $\M(v)$. For the sake of brevity we will
write $\mu$ for $\mu_\Gamma$, $i$ for $i_\Gamma$, and $\gaf$ for
$\otimes \gamma_i$. The sums below are always taken over $\GD \in G$.
\begin{equation*}
  \begin{split}
&\mu^* i^* st_* (ev^* \gaf \exp(\kaf \sf) \cap 
     [\M_{g,n}(V,\beta)]^{virt}) \\
&=\mu^*\st_* i^! (ev^* \gaf \exp(\kaf\sf) \cap
     [\M_{g,n}(V,\beta)]^{virt}) \\
&=\mu^* \st_* (i(\Gamma)^* (ev^* \gaf \exp(\kaf\sf)) \cap
     i^! [\M_{g,n}(V,\beta)]^{virt}) \\
&=\mu^*\st_* (i(\Gamma)^* (ev^* \gaf \exp(\kaf\sf)) \cap
     \frac{1}{N} \mu(\Gamma)_* \Delta_1^! \sum [X_\GD]^{virt}) \\
&=\frac{1}{N} \mu^* \st_* \mu(\Gamma)_* (\mu(\Gamma)^*
     i(\Gamma)^* (ev^* \gaf \exp(\kaf\sf)) \cap
     \Delta_1^! \sum [X_\GD]^{virt}) \\
&=\st_* \sum (\Delta^* (\otimes_{v\in \GD} ev^*\gaf_v 
     \exp(\kaf_v \sf) ) \cap \Delta_1^! \sum [X_\GD]^{virt}) \\
&=\sum \otimes_{v\in \GD} \st_* (ev^*\gaf_v 
     \exp(\kaf_v \sf) \cap \Delta_* \Delta_1^! [X_\GD]^{virt} ).  
  \end{split}
\end{equation*}
Summing over all $\beta$ gives the statement of the theorem taking
into account that $\Delta_* \Delta_1^!$ is the cap-product with the
Poincar\'e dual of the diagonal in $V^{E(\Gamma)} \times
V^{E(\Gamma)}$. 
\end{proof}

Th.~\ref{thm:cohft} provides a \cft\ determined by the $\kappa$
classes. One can similarly construct a \cft\ determined by the $\psi$
classes. Its potential is the usual potential. A more general
construction will appear in \cite{KK3}. 

\begin{rem}
  The potential of the \cft\ defined in the previous theorem coincides
  with the usual notion of potential of Gromov--Witten invariants up
  to terms quadratic in the variables $\bx$ which correspond to
  contributions from the moduli spaces $\M_{g,n}(V)$ where $2g-2+n\le
  0$. 
\end{rem}

\section{The Change of Coordinates}
\label{change}

\bigskip In this section we prove the change of coordinate formula on
the large phase space. Throughout the rest of this section, we fix a
homogeneous basis $\{ e_\alpha \}$ where $\alpha\in A$ of
$H^\bullet(V)$ such that $e_0$ is the identity element. We also fix a
total ordering on $A$.

\begin{rem}
  In this section we will not use the tautological classes
  $\kappa_{-1}(\gamma)$. 
\end{rem}

\begin{df}
  Let $\beta \in H_2^+(V,\nz)$, and $e_\alpha$, $\alpha=0, \ldots, r-1$
  be a basis of $H^\bullet(V)$. Assume that all $d_i >0$ and all $a_i
  \ge 0$. We define
\[
\begin{split}
&\la \sigma_{\nu_1} \ldots \sigma_{\nu_n} \,
\tau_{d_1,\mu_1} \ldots \tau_{d_k,\mu_k} \,
\kappa_{a_1,\alpha_1} \ldots \kappa_{a_l,\alpha_l} \ra_{g,\beta} := \\
&\int_{[\M_{g,n}(V,\beta)]^{virt}}
ev_1^*(e_{\nu_1}) \ldots ev_n^*(e_{\nu_n}) \\
&\hspace{3cm}
\pi_*(\psi_{n+1}^{d_1} ev_{n+1}^*(e_{\mu_1}) \ldots
\psi_{n+k}^{d_k} ev_{n+k}^*(e_{\mu_k}) )
\kappa_{a_1,\alpha_1} \ldots \kappa_{a_l,\alpha_l},
\end{split}
\]
where $\pi:\M_{g,n+k}(V,\beta) \to \M_{g,n}(V,\beta)$ ``forgets'' the last 
$k$ marked points. 
\end{df}

\begin{rem}
This definition differs from the standard one. However, if no $\kappa$ 
classes are present, then the intersection number above is the
standard intersection number of the $\psi$ and the pull-back classes
with $[\M_{g,n+k}(V,\beta)]^{virt}$. This definition is motivated by
the representation of the large phase space on the level of cohomology 
classes in Sec.~\ref{cohft}. Also, it will be easier to work with this
definition to derive the coordinate change below. 
\end{rem}

Let the sequence $\nu_1,\ldots, \nu_n$ contain $r_\nu$ elements $\nu$,
$\nu\in A$, the sequence $(d_{1},\mu_{1}), (d_{2}, \mu_{2}), \ldots,
(d_k, \mu_k)$ contain $m_{d,\mu}$ pairs $(d,\mu)$, where $d > 0$,
$\mu\in A$, and the sequence $(a_1, \alpha_1), (a_2, \alpha_2),
\ldots, (a_l, \alpha_l)$ contain $p_{a,\alpha}$ pairs $(a,\alpha)$
where $a\ge 0$, $\alpha\in A$. Then we also denote the intersection
number above by $\la \sif^\rf \tauf^\mf \kaf^\pf \ra_{g,\beta}$. One
has to be careful if $H^\bullet(V)$ has elements of odd degree. In
this case $\la \sif^\rf \tauf^\mf \kaf^\pf \ra_{g,\beta}$ denotes the
intersection number above with the following ordering. If $i<j$ then,
using the chosen order on $A$, 

a) $\nu_i \le \nu_j$; 

b) $d_i < d_j$, or $d_i=d_j$ and $\mu_i \le \mu_j$; 

c) $a_i < a_j$, or $a_i=a_j$ and $\alpha_i \le \alpha_j$.

\begin{df}
  We define
\[
\la \sif^\rf \tauf^\mf \kaf^\pf \ra_{g}  :=
\sum_{\beta\in H_2^+(V,\nz)} \la \sif^\rf \tauf^\mf \kaf^\pf
\ra_{g,\beta} \ q^\beta,  
\]
where $q$ is a formal variable. 
\end{df}

In the sequel we will consider the following collection of formal
variables: $\xf = (x^\nu)$, $\tf = (t_i^\mu)$, $\sf = (s_a^\alpha)$,
$i>0$, $a\ge 0$, $\nu,\mu,\alpha \in A$. These variables have the
following degrees: $|x^\nu|=|\nu|-2$, $|t_i^\mu|=2(i-1)+|\mu|$, and
$|s_a^\alpha|=2a+|\alpha|$. Note that the $\nz/2\nz$-degree is
determined by the upper index. Let
\[
\xf^\rf := \prod_{\nu} (x^\nu)^{r_\nu}, \quad
\tf^\mf := \prod_{d,\mu} (t_d^\mu)^{m_{d,\mu}}, \quad
\sf^\pf := \prod_{a,\alpha} (s_a^\alpha)^{p_{a,\alpha}}.
\]
Again one has to exercise care in case there are variables of odd
degree. In this case we order the products above so that 

$x^{\nu_1}$ preceeds $x^{\nu_2}$ if $\nu_1 \ge \nu_2$;

$t_{d_1}^{\mu_1}$ preceeds $t_{d_2}^{\mu_2}$ if $d_1>d_2$, or
$d_1=d_2$ and $\mu_1\ge \mu_2$; 

$s_{a_1}^{\alpha_1}$ preceeds $s_{a_2}^{\alpha_2}$ if $a_1>a_2$, or
$a_1=a_2$ and $\alpha_1 \ge \alpha_2$. \\
That is, we require the order on the products to be the opposite to the 
order on the intersection numbers. 

\begin{df}
  We define $\K_g \in \nc[[\Lambda,\xf,\tf,\sf]]$ by 
\[
\K_g(\xf,\tf,\sf) := \sum_{\rf,\mf,\pf}
\la \sif^\rf \tauf^\mf \kaf^\pf \ra_{g}
\, \frac{\sf^\pf}{\pf !} 
\, \frac{\tf^\mf}{\mf !}
\, \frac{\xf^\rf}{\rf !}
\]
\end{df}
where $\pf\,:=\,\prod_\nu\,p_\nu$ and $\pf !\,:=\,\prod_\nu\,p_\nu !$ (and
similarly for $\mf$ and $\rf$.

\begin{rem}
  In the above definition one could have chosen an arbitrary ordering
  for the intersection numbers, and then chosen the opposite ordering
  on the corresponding variables.
\end{rem}

The various degrees chosen for the variables together with the dimensions of
the cohomology classes and the virtual fundamental class insures the $\K_g$
has degree $2 (3-d) (1-g)$.

Note that $\K(\xf,\tf,\zef) = \F(\xf,\tf)$, the standard large phase
space potential if one sets $x^\nu = t_0^\nu$. Similarly,
$\K(\xf,\zef,\sf) = \G(\xf,\sf)$, the potential of the family of
\cfts\ determined by the $\kappa$ classes including the terms with
$2g-2+n\le 0$.

\begin{thm} \label{thm:change}
  Let $\tf(\sf)$ be determined by the following equation in $H^\bullet (V)$: 
\begin{equation} \label{eq:change}
e_0 - \sum_{d\ge 1} \, \theta^{d-1} t_d^\mu e_\mu =
\exp \left(- \sum_{a\ge 0} \theta^a s_a^\alpha e_\alpha \right), 
\end{equation}
where $\theta$ is an even formal parameter. Then $\F_g(\xf,\tf(\sf)) =
\G_g(\xf,\sf)$ for every $g\ge 0$.
\end{thm}

\begin{rem}
  In case when $V=pt$ Thm.~\ref{thm:change} reduces to Thm.~4.1 from
  \cite{MaZo}. The polynomials $t_a(\sf)$ are the Schur polynomials.
\end{rem}

We will prove the above theorem in a sequence of lemmata.

\begin{lm}
  Let $I$ and $J$ be two sets such that $I\cap J = \{ 1,\ldots,n\}$
  and $I\cup J = \{ 1,\dots, n+N \}$. Let $I':=I-\{ 1,\ldots,n \}$ and
  $J':=J-\{ 1,\ldots,n\}$. Consider the following commutative diagram
\[
\begin{CD}
  \M_{g,I}(V,\beta) @<\rho<< \M_{g,n+N}(V,\beta) \\
  @VV\pi V  @VV\pi V \\
  \M_{g,n}(V,\beta)   @<\rho<< \M_{g,J}(V,\beta), \\
\end{CD}
\]
  where the horizontal morphisms $\rho$ ``forget'' the marked points from 
  $J'$, and the vertical morphisms $\pi$ ``forget'' the marked points
  from $I'$. The morphisms $\pi$ and $\rho$ are flat, and
\[
\rho^*\pi_* = \pi_*\rho^* \quad \text{and} \quad 
\rho_* \pi^* = \pi^* \rho_*
\]
\end{lm}

\begin{proof}
  The morphisms $\pi$ and $\rho$ are flat as compositions of flat
  morphisms. Let $I$ be $\{1,\ldots,n+1\}$, and $J$ be $\{
  1,\ldots,n,n+2 \}$.  Then the commutative diagram above is close to
  a fibered square in the sense of Def.~\ref{df:fibered}
  (cf.~\cite{AC}).  Therefore one has $\rho^*\pi_* = \pi_*\rho^*$
  and $\rho_* \pi^* = \pi^* \rho_*$. Iterating, one obtains the
  statement of the lemma. 
\end{proof}

Consider the universal curve $\pi:\M_{g,n+1}(V,\beta) \to
\M_{g,n}(V,\beta)$. It has $n$ canonical sections $\sigma_i$, and each
of these sections is a regular embedding of codimension one.
Therefore, the image of each of these sections determines a Cartier
divisor on $\M_{g,n+1}(V,\beta)$. We denote the corresponding Chern
classes by $D_{i,n+1}\,\in\, H^2(\M_{g,n+1}(V,\beta))$. Equivalently,
$D_{i,n+1} = \sigma_{i*} 1$. The equalities below hold in
$H^\bullet(\M_{g,n+1}(V,\beta))$:
\[
\begin{split}
  & D_{i,n+1} D_{j,n+1} = 0 \quad \text{if} \ i\ne j,\\
  & \psi_i D_{i,n+1} = \psi_{n+1} D_{i,n+1} = 0.
\end{split}
\]
In addition, $\sigma_i^* D_{i,n+1} = -\psi_i$. 

Let $\pi:\M_{g,n+1}(V,\beta) \to \M_{g,n}(V,\beta)$ be the universal
curve. In the next two lemmas we will use the following properties.
Firstly, $\pi^* \psi_i = \psi_i - D_{i,n+1}$ proved in
\cite[Prop.~11]{Ge}. It follows that $\pi^* \psi_i^a = \psi_i^a +
(-1)^a D_{i,n+1}^a$. Secondly, $\pi^* ev_i^* = ev_i^*$.

The following lemma and its proof are similar to those in
\cite[Sec.~1]{AC}.

\begin{lm} \label{lm:kappalift}
  If $\gamma\in H^\bullet(V)$, then $\pi^* \kappa_a(\gamma) =
  \kappa_a(\gamma) - \psi_{n+1}^a ev_{n+1}^* \gamma$. 
\end{lm}

\begin{proof}
  Consider the following commutative diagram close to a fibered
  square: 
\[
\begin{CD}
  \M_{g,n+1}(V,\beta) @<\rho<< \M_{g,n+2}(V,\beta) \\
  @VV\pi V  @VV\pi V \\
  \M_{g,n}(V,\beta)   @<\rho<< \M_{g,J}(V,\beta), \\
\end{CD}
\]
where $J=\{ 1,\ldots,n,n+2\}$. Let $\sigma:\M_{g,n+1}(V,\beta) \to
\M_{g,n+2}(V,\beta)$ associated to the $(n+1)^{\text{st}}$ marked
point.  One has
\[
\begin{split}
& \pi^* \kappa_a(\gamma) 
  = \pi^*\rho_* (\psi_{n+2}^{a+1} ev_{n+2}^* \gamma) 
  = \rho_* \pi^* (\psi_{n+2}^{a+1} ev_{n+2}^* \gamma) \\
=& \rho_* (\psi_{n+2}^{a+1} ev_{n+2}^* \gamma) 
  + (-1)^{a+1} \rho_* (D_{n+1,n+2,}^{a+1} ev_{n+2}^* \gamma) \\
=& \kappa_a (\gamma) 
  + (-1)^{a+1} \rho_* \sigma_* 
  \sigma^*(D_{n+1,n+2,}^{a} ev_{n+2}^* \gamma) 
  = \kappa_a (\gamma) - \psi_{n+1}^a ev_{n+1}^* \gamma. 
\end{split}
\]

\end{proof}

\begin{df}
  Let $\gamma\in H^\bullet(V)$. Define the homomorphism
  $\varkappa_a(\gamma): H^\bullet(\M_{g,n+1}(V,\beta)) \to
  H^\bullet(\M_{g,n}(V,\beta))$ by
\[
\varkappa_a(\gamma) (x) := 
\pi_* (\psi_{n+1}^{a+1} ev_{n+1}^* \gamma \ x). 
\]
\end{df}

Note that $\varkappa_a (\gamma) (1) = \kappa_a(\gamma)$. 

\begin{lm} \label{lm:iterate}
  Assume that all $a_i > 0$ for $i=1,\ldots,N$. Then
\[
\varkappa_{a_1-1}(\gamma_1) \ldots \varkappa_{a_N-1}(\gamma_N) (x) = 
\pi_{N*} (\psi_{n+1}^{a_1} ev_{n+1}^* \gamma_1 \ldots
\psi_{n+N}^{a_N} ev_{n+N}^* \gamma_N \ x),
\]
where $\pi_N: \M_{g,n+N}(V,\beta) \to \M_{g,n}(V,\beta)$ ``forgets''
the last $N$ marked points. 
\end{lm}

\begin{proof}
  We proceed by induction. When $N=1$ the statement of the lemma is
  trivial. Assume that the statement is true for $N$ and prove it for
  $N+1$. We denote by $\pi'$ the universal curve
  $\M_{g,n+N+1}(V,\beta) \to \M_{g,n+N}(V,\beta)$. One has
\[
\begin{split}
&  \varkappa_{a_1-1}(\gamma_1) \ldots
   \varkappa_{a_{N+1}-1}(\gamma_{N+1})    (x) \\
=& \varkappa_{a_1-1}(\gamma_1) \ldots \varkappa_{a_N-1}(\gamma_N)
   (\pi'_* (\psi_{n+N+1}^{a_{N+1}} ev_{n+N+1}^* \gamma_{N+1}\ x)) \\
=& \pi_{N*} (\psi_{n+1}^{a_1} ev_{n+1}^* \gamma_1 \ldots
   \psi_{n+N}^{a_N} ev_{n+N}^* \gamma_N \ \pi'_*
   (\psi_{n+N+1}^{a_{N+1}} ev_{n+N+1}^* \gamma_{N+1} \ x)) \\
=& \pi_{N+1*} ( (\psi_{n+1}-D_{n+1,n+N+1})^{a_1} ev_{n+1}^* \gamma_1 
   \ldots \\
& \hspace{1cm} 
   (\psi_{n+N}-D_{n+N,n+N+1})^{a_N} ev_{n+N}^* \gamma_N \
   \psi_{n+N+1}^{a_{N+1}} ev_{n+N+1}^* \gamma_{N+1} \ x) ) \\
=& \pi_{N+1*} (\psi_{n+1}^{a_1} ev_{n+1}^* \gamma_1 \ldots
   \psi_{n+N}^{a_N} ev_{n+N}^* \gamma_N \ 
   \psi_{n+N+1}^{a_{N+1}} ev_{n+N+1}^* \gamma_{N+1} \ x). 
\end{split}
\]
\end{proof}

It follows from Lem.~\ref{lm:iterate} that the operators
$\varkappa_a(\gamma)$ super-commute. 

\begin{rem}
  If some if the numbers $a_i$ are equal to zero, then
  Lem.~\ref{lm:iterate} does not necessarily hold.
\end{rem}

\begin{lm} \label{lm:rel}
  Assume the conditions of Lem.~\ref{lm:iterate}. Then
\[
\begin{split}
&\pi_{N*} (\psi_{n+1}^{a_1} ev_{n+1}^* \gamma_1 \ldots
          \psi_{n+N}^{a_N} ev_{n+N}^* \gamma_N) \kappa_a (\gamma) \\
=&  \pi_{N+1*} (\psi_{n+1}^{a_1} ev_{n+1}^* \gamma_1 \ldots
             \psi_{n+N}^{a_N} ev_{n+N}^* \gamma_N \ 
             \psi_{n+1+N}^{a+1} ev_{n+N+1}^* \gamma) \\
-& \sum_{i=1}^k (-1)^{|\gamma_{i+1} \ldots \gamma_N||\gamma|} \pi_{N*} 
    (\psi_{n+1}^{a_1} ev_{n+1}^* \gamma_1 \ldots
     \psi_{n+i}^{a_i+a} ev_{n+i}^* (\gamma_i \gamma) \ldots
     \psi_{n+N}^{a_N} ev_{n+N}^* \gamma_N). 
\end{split}
\]
\end{lm}

\begin{proof}
  We proceed by induction. Let $N=1$, and $\pi_N = \pi$. Then, using
  Lem.~\ref{lm:kappalift} and Lem.~\ref{lm:iterate}, one gets
\[
  \pi_* (\psi_{n+1}^{a_1} ev_{n+1}^* \gamma_1 \pi^* \kappa_a (\gamma)) 
= \varkappa_{a_1-1}(\gamma_1) (\kappa_a (\gamma))
- \pi_* (\psi_{n+1}^{a_1+a} ev_{n+1}^* (\gamma_1 \gamma)). 
\]
This proves the statement of the lemma when $N=1$.  Now assume that
the statement is true for $N'=N-1$ and prove it for $N$. Denote by
$\pi'_N$ the natural morphism $\M_{g,n+N}(V,\beta) \to
\M_{g,n+1}(V,\beta)$.
\[
\begin{split}
&  \pi_{N*} (\psi_{n+1}^{a_1} ev_{n+1}^* \gamma_1 \ldots
          \psi_{n+N}^{a_{N}} ev_{n+N}^* \gamma_{N}) 
          \kappa_a (\gamma) \\
=& \pi_* (\psi_{n+1}^{a_1} ev_{n+1}^* \gamma_1 
    \pi'_{N*} (\psi_{n+2}^{a_2} ev_{n+2}^* \gamma_2 \ldots
              \psi_{n+N}^{a_{N}} ev_{n+N}^* \gamma_{N}) \\
&\hspace{3cm} (\kappa_a (\gamma) - \psi_{n+1}^a ev_{n+1}^* \gamma)). 
\end{split}
\]
The rest follows applying the induction hypothesis to the product 
\[
\pi'_{N*} (\psi_{n+2}^{a_2} ev_{n+2}^* \gamma_2 \ldots
\psi_{n+N}^{a_{N}} ev_{n+N}^* \gamma_{N}) \kappa_a (\gamma) ,
\]
and using Lem.~\ref{lm:iterate}. 
\end{proof}

Lemma \ref{lm:rel} provides a recursion relation for the intersection
numbers of the $\psi$ and the $\kappa$ classes. Let $\{ e_\alpha \}$,
$\alpha=0,\ldots,r$, be the chosen basis of $H^\bullet(V)$. Define
$c_{\alpha_1,\ldots,\alpha_j}^\mu$ by the formula
\[
e_{\alpha_1} \ldots e_{\alpha_j} = 
c_{\alpha_1,\ldots,\alpha_j}^\mu e_\mu. 
\]
(We assume summation over the repeating indices.) In particular,
$c_\alpha^\mu = \delta_\alpha^\mu$.

The following recursion relation follows from Lem.~\ref{lm:rel}. 
\[
\begin{split}
&\la \sif^\rf \tau_{d_1,\mu_1} \ldots \tau_{d_k,\mu_k} 
     \kappa_{a,\alpha} \kaf^\pf \ra_{g,\beta} 
= \la \sif^\rf \tau_{d_1,\mu_1} \ldots \tau_{d_k,\mu_k} 
     \tau_{a+1,\alpha} \kaf^\pf \ra_{g,\beta} \\
& \hspace{1cm} -\sum_{i=1}^k 
  (-1)^{|e_{\mu_{i+1}}\ldots e_{\mu_k}||e_\alpha|} 
  c_{\mu_i,\alpha}^\mu
  \la \sif^\rf \tau_{d_1,\mu_1} \ldots 
  \tau_{d_i+a,\mu} \ldots
  \tau_{d_k,\mu_k} \kaf^\pf \ra_{g,\beta}. 
\end{split}
\]
Note that the above relations also holds if one replaces $\la \ldots
\ra_{g,\beta}$ with $\la \ldots \ra_{g}$. 

It turns out that the equation above implies that for each $g\ge 0$
\begin{equation} \label{eq:st}
\frac{\del \K_g}{\del s_a^\alpha} =
\frac{\del \K_g}{\del t_{a+1}^\alpha} - 
\sum_{i=1}^\infty c_{\nu,\alpha}^\mu  t_i^\nu
\frac{\del \K_g}{\del t_{i+a}^\mu}. 
\end{equation}
We leave to the reader to check that all signs agree. 

Let us introduce the following standard notation: 
\[
\tw_i^\nu := 
\begin{cases}
  t_i^\nu & \text{unless $i=1$ and $\nu=0$}, \\
  t_1^0-1 & \text{if $i=1$ and $\nu=0$}. 
\end{cases}
\]
Then one can rewrite \eqref{eq:st} as 
\begin{equation} \label{eq:main}
\frac{\del \K_g}{\del s_a^\alpha} = 
-\sum_{i=1}^\infty c_{\nu,\alpha}^\mu  \tw_i^\nu
\frac{\del \K_g}{\del t_{i+a}^\mu}. 
\end{equation}

\begin{proof}
  (of Thm.~\ref{thm:change}.) We assume that $\tf(\sf)$ is determined
  by \eqref{eq:change}. It follows that $\tf(\zef)=\zef$, and 
  \begin{equation*}
    -\sum_{d\ge 1} \theta^{d-1} 
        \frac{\del t_d^\mu}{\del s_a^\alpha} e_\mu 
    = - \theta^a e_\alpha 
        \exp \left(- \sum_{a_1\ge 0} \theta^{a_1} 
        s_{a_1}^{\alpha_1} e_{\alpha_1} \right)
    = \sum_{d\ge a+1} \theta^{d-1} \tw_{d-a}^\mu e_\mu e_\alpha. 
  \end{equation*}
  It follows that for each $d$, $a$, and $\alpha$ such that $d\ge a+1$ 
  one has
\begin{equation} \label{eq:term}
  \frac{\del t_d^\mu}{\del s_a^\alpha} e_\mu 
  = - \tw_{d-a}^\nu e_\nu e_\alpha
  \quad \text{or, equivalently,} \quad
  \frac{\del t_d^\mu}{\del s_a^\alpha} 
  = - \tw_{d-a}^\nu c_{\nu,\alpha}^\mu, 
\end{equation}
and $\del t_d^\mu / \del s_a^\alpha = 0$ if $d \le a$. 

Consider the function $\K_g(\xf,\tf(\sf_0 + \sf), -\sf)$, where
$\sf_0$ is a constant. Differentiating it with respect to $s_a^\alpha$
provides using \eqref{eq:term} and \eqref{eq:main}
\begin{equation*}
\begin{split}
  &\frac{\del}{\del s_a^\alpha} 
  \left[ \K_g(\xf,\tf(\sf_0 + \sf), -\sf) \right] \\
  & = - \frac{\del \K_g}{\del s_a^\alpha} 
     (\xf,\tf(\sf_0 + \sf), -\sf) 
   + \sum_{d\ge 1} 
      \frac{\del t_d^\mu}{\del s_a^\alpha} (\sf_0 + \sf)
      \frac{\del \K_g}{\del t_d^\mu}
     (\xf,\tf(\sf_0 + \sf), -\sf)   \\
  & = - \frac{\del \K_g}{\del s_a^\alpha} 
     (\xf,\tf(\sf_0 + \sf), -\sf) 
   - \sum_{d\ge 1} 
     c_{\nu,\alpha}^\mu \tw_{d}^\nu (\sf_0 + \sf) 
     \frac{\del \K_g}{\del t_{d+a}^\mu}
     (\xf,\tf(\sf_0 + \sf), -\sf) 
    = 0. \\
\end{split} 
\end{equation*}
Therefore $\K_g(\xf,\tf(\sf_0 + \sf), -\sf)$ does not depend on
$\sf$. It follows that for all values $\sf_1$ and $\sf_2$ one has
\[
\K_g(\xf,\tf(\sf_1 + \sf_2), \zef) = 
\K_g(\xf,\tf(\sf_1), \sf_2) = 
\K_g(\xf,\zef, \sf_1+ \sf_2). 
\]
In particular, $\F_g(\xf,\tf(\sf)) = \G_g(\xf,\sf)$ for every $g\ge 0$.
\end{proof}

\begin{rem}
  Note that the condition $\tf(\sf)=\zef$ and \eqref{eq:term} are
  equivalent to \eqref{eq:change}, and determine $\tf(\sf)$
  completely. Note also that the coordinate change given by
  \eqref{eq:change} is invertible. 
\end{rem}

\begin{rem}
  The function $\tf(\sf)$ has the following Taylor coefficients: 
\[
\frac{\del^k t_d^\mu}{\del s_{a_1}^{\alpha_1} \ldots 
     \del s_{a_k}^{\alpha_k} }\mid_{\sf=\zef} = 
(-1)^{k+1} c_{\alpha_1, \ldots, \alpha_k}^\mu 
    \delta_{d,a_1+\ldots + a_k +1}. 
\]
\end{rem}

\section{Topological Recursion Relations}
\label{trr}

In this section we will derive the topological recursion relations for
$\G_0$ and $\G_1$ using the change of coordinates formula
\eqref{eq:change}. In \cite{KK} we represented the cohomology classes
by graphs to obtain the topological recursion relations when $V$ is a
convex variety and genus $g=0$. However, we were not able to extend
this technique to the general case since it is not clear that one can
pull back homology classes w.r.t. $\mu(\Gamma)$ from
Sec.~\ref{define}.

We will use the fact that $\G_g (\xf,\sf) = \F_g (\xf,\tf(\sf))$, where
$\tf(\sf)$ is determined by \eqref{eq:change}. Notice that $\del \G_g
/ \del x^\alpha = \del \F_g / \del x^\alpha$ since the coordinate
change $\tf(\sf)$ does not depend on $\xf$. We will raise and lower indices
in the usual manner.

\begin{prop}
  Let $a\ge 1$. Then
\[
\frac{\del^3 \G_0}{\del s_a^\alpha \del x^\mu \del x^\nu} = 
\frac{\del^2 \G_0}{\del s_{a-1}^\alpha \del x^\rho} \ 
\frac{\del^3 \G_0}{\del x_\rho \del x^\mu \del x^\nu}. 
\]
\end{prop}

\begin{proof}
Applying the chain rule one gets: 
  \begin{equation*}
\frac{\del^3 \G_0}{\del s_a^\alpha \del x^\mu \del x^\nu} =
\sum_{d\ge a+1} \frac{\del t_d^\xi}{\del s_a^\alpha} \ 
   \frac{\del^3 \F_0}{\del t_d^\xi \del x^\mu \del x^\nu} = 
\sum_{d\ge a+1} \frac{\del t_d^\xi}{\del s_a^\alpha} \ 
\frac{\del^2 \F_0}{\del t_{d-1}^\xi \del x^\rho}\
\frac{\del^3 \F_0}{\del x_\rho \del x^\mu \del x^\nu}. 
  \end{equation*}
The second equation uses that $\F_0$ satisfies the topological recursion 
relations. Similarly, 
\begin{equation*}
\frac{\del^2 \G_0}{\del s_{a-1}^\alpha \del x^\rho} \ 
\frac{\del^3 \G_0}{\del x_\rho \del x^\mu \del x^\nu} =
\sum_{d\ge a} \frac{\del t_d^\xi}{\del s_{a-1}^\alpha} \ 
\frac{\del^2 \F_0}{\del t_{d}^\xi \del x^\rho}\
\frac{\del^3 \F_0}{\del x_\rho \del x^\mu \del x^\nu}. 
\end{equation*}
Equation \eqref{eq:term} implies that $\del t_d^\xi / \del s_{a-1}^\alpha
= \del t_{d+1}^\xi / \del s_a^\alpha$, and the proposition follows. 
\end{proof}

\begin{prop} Let $|\alpha| \le 2$. Then
\begin{equation} \label{eq:deg2}
\frac{\del^3 \G_0}{\del s_0^\alpha \del x^\mu \del x^\nu} = 
\mathbf{D}_\alpha (\frac{\del \G_0}{\del x^\rho})
\frac{\del^3 \G_0}{\del x_\rho \del x^\mu \del x^\nu} +
c_{\rho,\alpha}^\xi x_\xi 
\frac{\del^3 \G_0}{\del x_\rho \del x^\mu \del x^\nu}
\end{equation}
where the differential operator $\mathbf{D}_\alpha$ is the
$\nc[[\bx,\bt,\bs]]$-linear operator defined by
\begin{equation}
\mathbf{D}_\alpha\,q^\beta\,:=\,q^\beta\,\int_\beta\,e_\alpha
\end{equation}
for all $\alpha$.
\end{prop}

\begin{proof}
We use the chain rule, \eqref{eq:term}, and the topological recursion
relations for $\F_0$: 
\begin{equation*}
\begin{split}
&\frac{\del^3 \G_0}{\del s_0^\alpha \del x^\mu \del x^\nu} 
= - \sum_{d\ge 1} c_{\xi,\alpha}^\zeta \tw_d^\xi 
\frac{\del^2 \F_0}{\del t_{d-1}^\zeta \del x^\rho}
\frac{\del^3 \F_0}{\del x_\rho \del x^\mu \del x^\nu}  \\
= & \mathbf{D}_\alpha (\frac{\del \F_0}{\del x^\rho})
\frac{\del^3 \F_0}{\del x_\rho \del x^\mu \del x^\nu} +
c_{\rho,\alpha}^\xi x_\xi 
\frac{\del^3 \F_0}{\del x_\rho \del x^\mu \del x^\nu}. \\
\end{split}
\end{equation*}
In the second equation we used the divisor equation for $\F$
\cite[2.6]{Ge2}. 
\end{proof}

\begin{rem}
The first term of the right hand side in \eqref{eq:deg2} contributes
only when $|e_\alpha|=2$. When $\alpha=0$ one can get \eqref{eq:deg2}
from 
\[
\frac{\del \G_0}{\del s_0^0} = \sum_{\rho} 
x^\rho \frac{\del \G_0}{\del x^\rho} - 2 \G_0. 
\]
This equation can be derived using the dilaton equation for $\F_0$
\cite[2.7]{Ge2}. 
\end{rem}

Similarly, one can derive the topological recursion relations for $\G$
in genus $1$ using known topological recursion relations for $\F$. We
state the results without proofs. 

\begin{prop}
  Let $a\ge 1$. Then 
  \begin{equation*}
\frac{\del \G_1}{\del s_a^\alpha} 
= \frac{\del^2 \G_0}{\del s_{a-1}^\alpha \del x^\rho}
\frac{\del \G_1}{\del x_\rho} 
+ \frac{1}{24} 
\frac{\del^3 \G_0}{\del s_{a-1}^\alpha \del x^\rho \del x_\rho}. 
  \end{equation*}
\end{prop}

\bibliographystyle{amsplain}

\providecommand{\bysame}{\leavevmode\hbox to3em{\hrulefill}\thinspace}

\end{document}